\pdfoutput=1
\documentclass[11pt]{extarticle}

\usepackage[left=1in,right=1in,top=1in,bottom=1.3in]{geometry}
\usepackage[utf8]{inputenc}
\usepackage[english]{babel}
\usepackage{lmodern}
\usepackage{amsmath,amssymb,amsthm,mathtools,indentfirst,xcolor,cite}

\usepackage[breaklinks,backref=page]{hyperref}
\usepackage{cleveref}
\hypersetup{
	colorlinks = true, 
	urlcolor = cyan, 
	linkcolor = teal, 
	citecolor = cyan 
}

\renewcommand*{\backref}[1]{}
\renewcommand*{\backrefalt}[4]{
	\ifcase #1 \textcolor{red}{Not cited.}%
	\or $\uparrow$#2%
	\else $\uparrow$#2%
	\fi%
}
\let\OLDthebibliography\thebibliography
\renewcommand\thebibliography[1]{
	\OLDthebibliography{#1}
	\setlength{\parskip}{0pt}
	\setlength{\itemsep}{0pt plus 0.3ex}
}

\newtheorem{Theorem}{Theorem}
\newtheorem{Corollary}{Corollary}
\newtheorem{Lemma}{Lemma}
\newtheorem{Claim}{Claim}
\crefname{Theorem}{theorem}{theorems}
\crefname{Corollary}{corollary}{corollaries}
\crefname{Lemma}{lemma}{lemmas}
\crefname{Claim}{claim}{claims}
\Crefname{Theorem}{Theorem}{Theorems}
\Crefname{Corollary}{Corollary}{Corollaries}
\Crefname{Lemma}{Lemma}{Lemmas}
\Crefname{Claim}{Claim}{Claims}

\newcommand{\R}{\mathbb R}
\newcommand{\Q}{\mathbb Q}
\newcommand{\N}{\mathbb N}
\newcommand{\MR}{\smash{\xrightarrow%
  [\smash{\raisebox{3pt}{\text{{\tiny \upshape MR}}}}]%
  {}}}

\newcommand{\norm}[1]{\lVert #1\rVert}
\newcommand{\tensq}[1]{\bigl(#1\bigr)^{[2]}}
\newcommand{\trans}{\mathsf T}

\begin{document}
\title{Ramsey sets are canonically Ramsey}
\author{Arsenii Sagdeev\thanks{Karlsruhe Institute of Technology, Karlsruhe, Germany, \href{mailto:sagdeevarsenii@gmail.com}{\tt sagdeevarsenii@gmail.com}.}}
\date{}
\maketitle

\begin{abstract}
A set \(A\) of points is called Ramsey if, for each \(r\), every
\(r\)-coloring of a sufficiently high-dimensional Euclidean space
contains a monochromatic congruent copy of \(A\).
We show that for each Ramsey set \(A\), every \(r\)-coloring of a
sufficiently high-dimensional Euclidean space, regardless of the number of colors used,
contains a monochromatic or a rainbow congruent copy of \(A\).
This answers a question raised in virtually every paper on the topic.
\end{abstract}

\noindent
{\textbf{AI disclosure.} This is an early announcement. Most of the proofs and the writing is due to ChatGPT-6 Astra, the introduction is due to the author. The exposition is in active preparation and will replace this draft shortly.}



\section{Introduction}\label{sec:intro}

In this paper, a \emph{copy} is a congruent copy.
For \(A \subset \R^d,\ X \subset \R^n\) and a positive integer $r$, 
we write \(X\xrightarrow{r} A\) if every \(r\)-coloring of \(X\)
contains a monochromatic copy of \(A\). A set $A \subset \R^d$ is called \textit{Ramsey}
if, for every positive integer \(r\), there is an \(n\) such that
\(\R^n\xrightarrow{r} A\).

Erd\H os, Graham, Montgomery,
Rothschild, Spencer and Straus introduced this notion and, among many other results, proved that every Ramsey set is finite~\cite[Theorem~19]{EGMRSS2} and spherical~\cite[Theorem~13]{EGMRSS}, i.e., lie on the surface of a hypersphere. Graham conjectured~\cite{G17} that these two necessary conditions are also sufficient, while Leader, Russell, and Walters~\cite{LRW} posed a `rival conjecture'.
Both were were recently refuted. First, P\'alv\"olgyi~\cite{DP26} constructed a $7$-point counterexample to the first one, while shortly after, OpenAI~\cite{OpenAI} found a $4$-point counterexample to the `rival conjecture'. Moreover, OpenAI~\cite{OpenAI} announced a complete characterization of Ramsey sets in terms of a tensor condition over their coordinate fields, see \Cref{thm:criterion} for a precise statement. We refer the reader to the work~\cite{ALS} of Axenovich, Liu and Sagdeev for the background on this subject and its rich bibliography. 

Each problem in Ramsey theory has a `canonical' counterpart where the number of colors is not fixed, see classical examples in~\cite{GST}. For \(A \subset \R^d,\ X \subset \R^n\),
we write \(X\MR A\) if every  coloring of \(X\), regardless of the number of colors used, contains a copy of $A$ that is either monochromatic or \textit{rainbow}, i.e., with points of pairwise different colors. We say that \(A \subset \R^d\) is
\emph{canonically Ramsey} if there is an \(n\) such that \(\R^n\MR A\).

Mao, Ozeki and Wang~\cite{MOW} introduced this class of problems and that $\R^3 \MR T$, where $T$ is a triangle with angles \(30^\circ\), \(60^\circ\), \(90^\circ\). (Here and in what follows, standard geometric shapes are viewed as their vertex sets.) Cheng and Xu~\cite{CX} as well as Geh\'er, Sagdeev and T\'oth~\cite{GST} strengthened this result by showing that $\R^3 \MR T$ if $T$ is a right or an acute triangle, respectively. Fang, Ge, Shu, Xu, Xu and Yang~\cite{Fang} proved that if $T$ is an obtuse triangle, then $\R^4 \MR T$. It is somewhat natural to expect the dimension in these results can be further reduced to 2.  Rauch-M\"uller and Sagdeev~\cite{RS} recently confirmed this for 36 specific triangles.
 
Cheng and Xu~\cite{CX} also proved that squares are canonically Ramsey. Their result was strengthened in two `independent' directions. First, Geh\'er, Sagdeev and T\'oth~\cite{GST} showed that hypercubes are canonically Ramsey. Second, Ge, Shu and Xu~\cite{GSX} proved that rectangles are canonically Ramsey. (The latter result was incorporated into~\cite{Fang}.) Later Shaw~\cite{ShawCuboids} found a common strengthening by showing that hyperrectangles are canonically Ramsey.

Ge, Shu, Xu and Yu~\cite{GSXY} proved that nondegenerate simplices are canonically Ramsey. Shortly afterwards, Shaw~\cite{Shaw} combined their ideas with his earlier argument for hyperrectangles to conclude that Cartesian products of nondegenerate simplices are canonically Ramsey. Finally, Shaw~\cite{ShawPolygons} proved that regular polygons
with a prime number of vertices and all their Cartesian powers are canonically Ramsey as well.

In fact, it seems hard to avoid rainbow copies of a given set $A$ if the number of colors used is much larger than the cardinality $|A|$. Therefore, it is natural to expect that Ramsey sets must also be canonically Ramsey. This intuition was first expressed in~\cite{GST}, appeared in~\cite{Fang} as a conjecture and was discussed in all subsequent `canonical' papers~\cite{GSXY, ShawCuboids, Shaw, ShawPolygons}. In the present paper, we confirm this intuition.

\begin{Theorem}\label{thm:main}
If sets \(A\) and \(B\) are Ramsey, then there is a positive integer \(n\) such that every coloring of \(\R^n\), regardless of the number of colors used, contains either a monochromatic copy of $A$ or a rainbow copy of $B$.
\end{Theorem}

\begin{Corollary}\label{cor:diagonal}
Every Ramsey set is canonically Ramsey.
\end{Corollary}

\noindent\textbf{Remark.} It seems that one can deduce from~\cite{OpenAI} that
\begin{equation} \label{eq:strong}
    \mbox{\it For every Ramsey set $A$ and every $r \in \N$, there is a Ramsey set $X$ such that $X \xrightarrow{r} A$.}
\end{equation}
Combined with our proof of \Cref{thm:main}, it would lead to a similar `Ramsey witness' strengthening of our result. Namely, if sets \(A\) and \(B\) are Ramsey, then there is another Ramsey set \(C\) such that every coloring of \(C\) contains either a monochromatic copy of $A$ or a rainbow copy of $B$. However, to the best of our knowledge, \eqref{eq:strong} has not yet been stated as a theorem anywhere, so we are not pursuing this strengthening in the present paper. We also remark that, since Ramsey sets are spherical, \eqref{eq:strong} would imply, in a strong form, Graham's conjecture that Ramsey sets are sphere-Ramsey, see Graham~\cite[Conjecture~11.3.5]{G17} and P\'alv\"olgyi~\cite{DP26ncs}.
\vspace{2.5mm}

\Cref{cor:diagonal} follows from \Cref{thm:main} by taking \(A=B\). However, \Cref{thm:main} also follows \Cref{cor:diagonal} by applying it to the Cartesian product $A \times B$. This relies on the classical result that Cartesian products of Ramsey sets are Ramsey, see~\cite[Theorem~20]{EGMRSS}.

We note that all previously known canonically Ramsey sets are Ramsey. From this perspective, \Cref{thm:main} is a common strengthening of all the `canonical' results mentioned above. Its proof is also a mix of earlier arguments, so in what follows, we try to attribute the ideas properly.

Finally, we announce that \textit{not all canonically Ramsey sets are Ramsey}. One particular example is $A=\{0,1,2\}\subset\R$. The question of whether this specific non-spherical (and thus non-Ramsey) set is canonically Ramsey was raised explicitly in virtually every paper on the topic. In Appendix~\ref{appA}, we present a sketch of the proof that this set is canonically Ramsey found by GPT-6 Astra. Currently we are working with Kenneth Moore on an alternative proof that will be published separately. The complete characterization of canonically Ramsey sets remains open.

In what follows, we write \([n]=\{1,\ldots,n\}\). Vectors are in bold,
\(u_i\) stands for the \(i\)-th coordinate of \(\mathbf{u}\) and similarly for other letters. We write $\mathbf{v}\cdot\mathbf{u}$ for the dot product of two vectors $\mathbf{v}, \mathbf{u}$ and $\norm{\mathbf{u}}$ for the Euclidean norm $\sqrt{\mathbf{u}\cdot\mathbf{u}}$. Cartesian products are naturally embedded into Euclidean spaces of larger dimension. We say that $A\subset \R^d$ is \textit{of full affine span} if no hyperplane in $\R^d$ contains $A$.

\vspace{2.5mm}
\noindent\textbf{Outline of the proof of \Cref{thm:main}.}
In \Cref{sec:geometry}, we give a `product-like' constriction containing many copies of a given set, all sharing a common point. \Cref{sec:criterion} uses OpenAI's characterization of Ramsey sets as a black box to show that the constructed set is Ramsey. Consider a coloring of a Cartesian product of sufficiently many copies of this construction, and suppose that no copy of \(B\) is rainbow. Every copy of \(B\) then has two points of the same color, which defines an auxiliary coloring of the product with $\binom{|B|}{2}$ colors. In \Cref{sec:line}, we apply the Hales--Jewett theorem to find a monochromatic line, which corresponds to a monochromatic copy of \(A\) at one of finitely many scales in the original coloring. In \Cref{sec:scale}, we repeat this argument sufficiently many times to find a monochromatic Cartesian power of a scaled copy of $A$, in which we will take a diagonal to find the desired monochromatic copy of $A$ at the original scale.

\section{Copies sharing a common point}\label{sec:geometry}

Informally speaking, every point $(\mathbf{v},\mathbf{u})$ of the Cartesian product $V \times U$ lies in the slice $V\times\{\mathbf{u}\}$ forming a copy of $V$ in $V \times U$. In this section, we give a more involved `product-like' construction with an additional property that all these copies share a common point. Similar construction was recently used by Moore~\cite{Moore}, who arranged congruent pyramids so that their apex vertices form a prescribed configuration. He used this construction to prove that a set obtained from a Ramsey set by adding another point outside its affine hull is also Ramsey, answering a question of Ivan, Leader and Walters~\cite{ILW}. Mirabi~\cite{Mirabi} obtained this result independently shortly afterwards. Our orthogonal height adjustment below was also used by Ge, Shu, Xu and Yu in their recent result that nondegenerate simplices are canonically Ramsey, see~\cite[Lemma~2.4]{GSXY}.

Let \(V\subset\R^n\) and \(U\subset\R^m\), and suppose that all points of \(U\) have the same norm \(r<1\). We view the space
\(\R^{n(m+1)}\) as \((\R\times\R^m)^n\), i.e.,  a point of \(\R^{n(m+1)}\)
has \(n\) blocks, each consisting of one real number and one vector
in \(\R^m\). Define \(\phi_r:\R^n\times \R^m\to \R^{n(m+1)}\) and $F(V,U) \subset \R^{n(m+1)}$ by
\begin{equation}\label{eq:family}
  \phi_r(\mathbf{v},\mathbf{u})
     =\bigl((\sqrt{1-r^2}\,v_i,\ v_i\mathbf{u})\bigr)_{i=1}^n,
  \qquad F(V,U)=\phi_r(V\times U).
\end{equation}
The following identities describe the `slices' of this construction.

\begin{Lemma}\label{lem:distances}
For \(\mathbf{v},\mathbf{v}'\in V\) and \(\mathbf{u},\mathbf{u}'\in U\),
\begin{align}
 \norm{\phi_r(\mathbf{v},\mathbf{u})-\phi_r(\mathbf{v}',\mathbf{u})}^2
    &=\norm{\mathbf{v}-\mathbf{v}'}^2, \label{eq:fixed-y}\\
 \norm{\phi_r(\mathbf{v},\mathbf{u})-\phi_r(\mathbf{v},\mathbf{u}')}^2
    &=\norm{\mathbf{v}}^2\,\norm{\mathbf{u}-\mathbf{u}'}^2. \label{eq:fixed-u}
\end{align}
\end{Lemma}

\begin{proof}
The \(i\)-th block of the difference $\phi_r(\mathbf{v},\mathbf{u})-\phi_r(\mathbf{v}',\mathbf{u})$ is
\((v_i-v_i')(\sqrt{1-r^2},\mathbf{u})\). Its squared norm is \((v_i-v_i')^2\) since \(1-r^2+\norm{\mathbf{u}}^2=1\). Summing over \(i\) gives~\eqref{eq:fixed-y}.

Similarly, the \(i\)-th block of the difference $\phi_r(\mathbf{v},\mathbf{u})-\phi_r(\mathbf{v},\mathbf{u}')$ is
\((0,v_i(\mathbf{u}-\mathbf{u}'))\). Its squared norm is
\(v_i^2\norm{\mathbf{u}-\mathbf{u}'}^2\). Summing over \(i\) gives~\eqref{eq:fixed-u}.
\end{proof}

In particular,  for every \(\mathbf{v} \in V\), the slice $\phi_r(\{\mathbf{v}\}\times U) \subset F(V,U)$ is a copy of \(U\) scaled by $\norm{\mathbf{v}}$. Similarly, for every \(\mathbf{u} \in U\), the slice $\phi_r(V\times\{\mathbf{u}\}) \subset F(V,U)$ is a copy of \(V\). In addition, if \(\mathbf{0}\in V\), then all these copies of $V$ share the point \(\phi_r(\mathbf{0},\mathbf{u})=\mathbf{0}\).

\section{Characterization of Ramsey sets}\label{sec:criterion}

In this section, we use OpenAI's characterization of Ramsey sets~\cite[Theorem~1.1]{OpenAI} as a black box to show the set $F(V,U)$ constructed in \Cref{sec:geometry} is Ramsey. We use an alternative form of the characterization that does not involve matrices with tensor entries, see~\cite[Proposition~2.1]{OpenAI} for their equivalence. However, this variant still involves tensor products. We recall their basic properties below and also refer the reader to Lang's classic textbook~\cite[Chapter~XVI]{Lang} for further algebraic background. Note that in the present paper, we use tensor products only in this section, while the other ones are purely combinatorial. 

We view $\R$ as an infinite-dimensional vector space over $\Q$. Consider the bilinear over $\Q$ function $\otimes_{\Q}$  that maps  two real numbers \(z,z'\) to the symbol \(z\otimes_{\Q}z'\). The tensor product $\R \otimes_{\Q} \R$ is then defined as the vector space over $\Q$ spanned by all these symbols with no additional relations but those following from the bilinearity.

Write \(z^{[2]}=z\otimes_{\Q}z\). To check an identity in
\(\R\otimes_{\Q}\R\) of the form \(\sum_i c_i z_i^{[2]}=0\) with finitely many \(z_i\in\R\)
and \(c_i\in\Q\), choose a basis \(\xi_1,\ldots,\xi_n\) of the
rational span of all \(z_i\), expand the product by bilinearity and check $n^2$ corresponding identities separately. Namely, if \(z_i=\sum_{j=1}^n
s_{ij}\xi_j\) with \(s_{ij}\in\Q\), then $\sum_i c_i z_i^{[2]}=0$ if and only if
\begin{equation}\label{eq:tensor-coordinates}
   \sum_i c_i s_{ij}s_{ik}=0
       \qquad \text{for all } 1\leq j,k \leq n.
\end{equation}
Indeed, expanding \(z_i\otimes_{\Q}z_i\) gives coefficient
\(s_{ij}s_{ik}\) in from of \(\xi_j\otimes_{\Q}\xi_k\),
and these basis tensors are linearly independent over \(\Q\) in
\(\R\otimes_{\Q}\R\). We will need the following two consequences of~\eqref{eq:tensor-coordinates}.

\begin{Lemma}\label{lem:tensor_con}
In the notation of the paragraph above, if \(\sum_i c_i z_i^{[2]}=0\), then
\begin{align}
 \sum_i c_i z_i^2 &=0, \label{eq:square1}\\
 \sum_i c_i\tensq{az_i}&=0 \quad \mbox{ for all } a \in \R. \label{eq:square2}
\end{align}
\end{Lemma}
\begin{proof}
To prove~\eqref{eq:square1}, multiply each identity
in~\eqref{eq:tensor-coordinates} by \(\xi_j\xi_k\) and sum
over \(j,k\).
For~\eqref{eq:square2}, assume that $a\neq 0$, since otherwise there is nothing to prove.  Now the numbers \(a\xi_1,\ldots,a\xi_n\) from a
basis of the rational span of all \(az_i\), with the same coefficients $s_{ij}$, and thus~\eqref{eq:tensor-coordinates} completes the proof.
\end{proof}

\begin{Theorem}[Characterization of Ramsey sets~{\cite{OpenAI}}]\label{thm:criterion}
A finite set \(U\subset\R^m\) of full affine span is Ramsey if and only if there are finitely many weights \(c_i\in\Q\)
affine functions
\(\ell_i:\R^m\to\R\), of the form
\(\ell_i(\mathbf{u})=b_i+\mathbf{w}_i\cdot\mathbf{u}\) with \(b_i\in\R\) and
\(\mathbf{w}_i\in\R^m\), such that
\begin{alignat}{2}
 \sum_i c_i\tensq{\ell_i(\mathbf{u})} &=0 \quad
       && \mbox{ for all }\mathbf{u}\in U, \label{eq:certificate}\\
 \sum_i c_i(\mathbf{w}_i\cdot\mathbf{x})^2 &=\norm{\mathbf{x}}^2
       \quad && \mbox{ for all }\mathbf{x}\in\R^m. \label{eq:certificate2}
\end{alignat}
\end{Theorem}

We call the functions $\ell_i$ and weights $c_i$ from the statement of \Cref{thm:criterion} a
\emph{certificate} for \(U\) and \(\mathbf{w}_i\cdot\mathbf{x}
=\ell_i(\mathbf{x})-\ell_i(\mathbf{0})\) the \emph{linear part} of \(\ell_i\)
evaluated at \(\mathbf{x}\). Note that \eqref{eq:certificate} is the
quadratic-form version of the condition that \(\sum_i c_i\mathbf{w}_i\mathbf{w}_i^{\trans}\) is an \(m\times m\) identity matrix.

\begin{Lemma}\label{lem:family}
If \(V\subset\R^n\) and \(U\subset\R^m\) are Ramsey
sets of full affine span, \(\mathbf{0}\in V\), and all points
of \(U\) have the same norm \(r<1\), then \(F(V,U)\) from~\eqref{eq:family}
is Ramsey.
\end{Lemma}

\begin{proof}
Let the functions $\ell_i(\mathbf{u})=b_i+\mathbf{w}_i\cdot\mathbf{u}$ and weights $c_i$ be a certificate for \(U\).  For \(\mathbf{u}\in U\), combine~\eqref{eq:square1} and~\eqref{eq:certificate} to get \(\sum_i c_i\ell_i(\mathbf{u})^2=0\).
Expand the squares and apply~\eqref{eq:certificate2} with \(\mathbf{x}=\mathbf{u}\) to get
\(\sum_i c_i b_i^2+2(\sum_i c_i b_i\mathbf{w}_i)\cdot\mathbf{u}+r^2=0\).
This affine function of \(\mathbf{u}\) vanishes on a set $U$ of full affine
span, so both its linear and constant terms vanish, i.e.,
\(\sum_i c_i b_i\mathbf{w}_i=\mathbf{0}\) and \(\sum_i c_i b_i^2=-r^2\).
For any \(a\in\R\) and \(\mathbf{z}\in\R^m\), expanding once more gives
\begin{equation}\label{eq:one-block}
\begin{aligned}
 \sum_i c_i(b_i a+\mathbf{w}_i\cdot\mathbf{z})^2
 &=a^2\sum_i c_i b_i^2
    +2a\Bigl(\sum_i c_i b_i\mathbf{w}_i\Bigr)\cdot\mathbf{z}
    +\sum_i c_i(\mathbf{w}_i\cdot\mathbf{z})^2\\
 &=-r^2a^2+\norm{\mathbf{z}}^2 .
\end{aligned}
\end{equation}

We now construct a certificate on \(\R^{n(m+1)}\). Write a general point
\(\mathbf{z}\in \R^{n(m+1)}\) as \(\bigl((\sqrt{1-r^2}\,y_j,\mathbf{x}_j)\bigr)_{j=1}^n\),
where \(\mathbf{y}=(y_1,\ldots,y_n)\in\R^n\) and \(\mathbf{x}_j\in\R^m\).
These coordinates are unique because \(r<1\). At a point of
\(F(V,U)\) the \(j\)-th block records \(v_j\) and \(v_j\mathbf{u}\).
We can therefore recover \(v_j\ell_i(\mathbf{u})=b_iv_j+\mathbf{w}_i\cdot(v_j\mathbf{u})\)
by a linear function of just this block. This gives the first group of
functions \(L_{ij}:\R^{n(m+1)}\to\R\), defined by
\(L_{ij}(\mathbf{z})=b_i y_j+\mathbf{w}_i\cdot\mathbf{x}_j\) for \(j\in[n]\),
with weights \(c_i\). Formula~\eqref{eq:one-block} accounts for the
\(\mathbf{x}_j\)-coordinates. To supply the remaining contribution in the
\(\mathbf{y}\)-coordinates, choose a certificate \((\lambda_{\kappa},\gamma_{\kappa})\) for \(V\),
where \(\lambda_{\kappa}:\R^n\to\R\) has the form
\(\lambda_{\kappa}(\mathbf{v})=\beta_{\kappa}+\boldsymbol{\omega}_{\kappa}\cdot\mathbf{v}\),
with \(\beta_{\kappa}\in\R\), \(\boldsymbol{\omega}_{\kappa}\in\R^n\), and
\(\gamma_{\kappa}\in\Q\). Add the affine functions
\(\Lambda_{\kappa}:\R^{n(m+1)}\to\R\) given by
\(\Lambda_{\kappa}(\mathbf{z})=\lambda_{\kappa}(\mathbf{y})
=\beta_{\kappa}+\boldsymbol{\omega}_{\kappa}\cdot\mathbf{y}\),
with weights \(\gamma_{\kappa}\).

We check~\eqref{eq:certificate} first. At a
point \(\mathbf{z}=\phi_r(\mathbf{v},\mathbf{u})\) of \(F(V,U)\) we have
\(\mathbf{y}=\mathbf{v}\) and \(\mathbf{x}_j=v_j\mathbf{u}\). Hence
\(L_{ij}(\mathbf{z})=v_j\ell_i(\mathbf{u})\).
For each fixed \(j\), the tensor sum over \(i\) vanishes by~\eqref{eq:square2}, applied
with \(a=v_j\). The tensor sum for the \(\Lambda_{\kappa}\) vanishes because
\(\Lambda_{\kappa}(\mathbf{z})=\lambda_{\kappa}(\mathbf{v})\) and \(\mathbf{v}\in V\).

For~\eqref{eq:certificate2}, take an arbitrary direction \(\mathbf{z}\in \R^{n(m+1)}\) in
the same block coordinates. The \(L_{ij}\) are linear, so their linear
parts are their values. By~\eqref{eq:one-block}, their contribution
in block \(j\) is \(-r^2y_j^2+\norm{\mathbf{x}_j}^2\). The linear parts
of the \(\Lambda_{\kappa}\) are \(\boldsymbol{\omega}_{\kappa}\cdot\mathbf{y}\), whose weighted squares
sum to \(\norm{\mathbf{y}}^2\) by the certificate for \(V\). The total is
\[
 \sum_{j=1}^n\bigl(-r^2y_j^2+\norm{\mathbf{x}_j}^2\bigr)+\norm{\mathbf{y}}^2
 =\sum_{j=1}^n\bigl((1-r^2)y_j^2+\norm{\mathbf{x}_j}^2\bigr)
 =\norm{\mathbf{z}}^2 .
\]

It remains to check the affine span. The origin belongs to \(F(V,U)\).
Suppose a linear function \(H:\R^{n(m+1)}\to\R\) vanishes on \(F(V,U)\). In the
above coordinates \(H(\mathbf{z})=\sum_j (a_j y_j+\mathbf{t}_j\cdot\mathbf{x}_j)\),
with \(a_j\in\R\) and \(\mathbf{t}_j\in\R^m\). Substitution of~\eqref{eq:family}
gives \(\sum_j v_j(a_j+\mathbf{t}_j\cdot\mathbf{u})=0\) for every \(\mathbf{v}\in V\)
and \(\mathbf{u}\in U\). For fixed \(\mathbf{u}\), the vectors \(\mathbf{v}\) span
\(\R^n\), so \(a_j+\mathbf{t}_j\cdot\mathbf{u}=0\) for each \(j\).
Since \(U\) affinely spans \(\R^m\), every \(a_j\) and \(\mathbf{t}_j\)
vanishes. Thus \(F(V,U)\) spans \(\R^{n(m+1)}\), and \Cref{thm:criterion}
shows that \(F(V,U)\) is Ramsey.
\end{proof}

\section{A monochromatic copy at one of finitely many scales}\label{sec:line}

For a positive integer \(h\), write \(A_h=h^{-1/2}A\) for \(A\) scaled
by \(1/\sqrt h\). We first force a monochromatic copy at one of finitely
many such scales whenever a coloring has no rainbow \(B\). For the
product argument in the last section, the set \(Z\) in which we force
this copy must itself be Ramsey.

\begin{Lemma}\label{lem:scales}
Let \(A\) and \(B\) be Ramsey configurations with at least two points
each. There are a positive integer \(q\) and a Ramsey configuration
\(Z\) such that every coloring of \(Z\) contains a rainbow copy of
\(B\) or a monochromatic copy of some \(A_h\), \(h\in[q]\).
\end{Lemma}

\begin{proof}
Work in the affine spans of \(B=\{\mathbf{b}_1,\ldots,\mathbf{b}_s\}\subset\R^d\)
and \(A\subset\R^f\), and use the block decomposition \(\R^{d(f+1)}=(\R\times\R^f)^d\) as
in~\eqref{eq:family}. Since every Ramsey configuration is
spherical~\cite[Theorem~13]{EGMRSS}, we may take the circumcenter of
\(A\) to be the origin. Denote its circumradius by \(\rho\).

Let \(P=\binom{[s]}2\). Think of a letter \(\theta\in A^P\) as a
function \(\theta:P\to A\) that assigns a point of \(A\) to each
pair of vertices of \(B\). By the Hales--Jewett theorem~\cite{HJ},
there is a positive integer \(N\) such that every \(|P|\)-coloring of
the words \((A^P)^N\) contains a monochromatic combinatorial line.
Geh\'er, Sagdeev and T\'oth~\cite[\S3]{GST} use the canonical
Hales--Jewett theorem for the maximum norm, where repeating a
coordinate preserves distances. In the Euclidean metric each repetition
contributes another summand to the squared distance.
A combinatorial line is obtained by choosing a nonempty set \(I\subseteq[N]\),
fixing the letters outside \(I\), and writing the same letter
\(\theta\) in every position of \(I\), once for each \(\theta\in A^P\).

\smallskip\noindent\emph{A copy for each word.}
First prepare a configuration for each \(p'=\{i',j'\}\in P\), with
\(i'<j'\). Choose a single positive integer \(q\), divisible by
\(1,\ldots,N\), large enough that \(\delta_{p'}\rho<1\) for every
\(p'\), where
\(\delta_{p'}=\sqrt{N|P|/q}\big/\norm{\mathbf{b}_{j'}-\mathbf{b}_{i'}}\).
All these distances are positive, so a sufficiently large multiple
of \(N!\) will do. The number of variable positions on a Hales--Jewett
line is unknown at this point. This scale problem appears in the work of
Leader, Russell and Walters~\cite[\S\S2--3]{LRW}. OpenAI~\cite[\S5.3]{OpenAI}
prepares for every possible number of variable positions before taking
the line. Here divisibility by \(1,\ldots,N\) will let us recover the
original scale for every such number.
By \Cref{lem:family},
\(F_{p'}=F(B-\mathbf{b}_{i'},\delta_{p'}A)\subset \R^{d(f+1)}\) is Ramsey.
Indeed, translations and positive dilations preserve the Ramsey
property, and \(B-\mathbf{b}_{i'}\) contains the origin.

Define \(\eta_{p'}:B\times A\to F_{p'}\) by
\(\eta_{p'}(\mathbf{b},\mathbf{a})=
\phi_{\delta_{p'}\rho}(\mathbf{b}-\mathbf{b}_{i'},\delta_{p'}\mathbf{a})\).
With \(\mathbf{a}\) fixed, this is an isometry in \(\mathbf{b}\)
by~\eqref{eq:fixed-y}. As \(\mathbf{a}\) varies, the image of
\(\mathbf{b}_{i'}\) is always the origin. For \(\mathbf{a},\mathbf{a}'\in A\),
the other endpoint satisfies
\begin{equation}\label{eq:moving-vertex}
\begin{aligned}
 \norm{\eta_{p'}(\mathbf{b}_{j'},\mathbf{a})-\eta_{p'}(\mathbf{b}_{j'},\mathbf{a}')}^2
 &=\norm{\mathbf{b}_{j'}-\mathbf{b}_{i'}}^2\,\delta_{p'}^2
       \norm{\mathbf{a}-\mathbf{a}'}^2\\
 &=\frac{N|P|}{q}\norm{\mathbf{a}-\mathbf{a}'}^2 .
\end{aligned}
\end{equation}
The first equality is~\eqref{eq:fixed-u}. The second explains the
choice of \(\delta_{p'}\). Every pair now has the same scale.

Set \(Z=(N|P|)^{-1/2}\bigl(\prod_{p'\in P} F_{p'}\bigr)^N\).
Its ambient space is \(\R^{N|P|d(f+1)}\),
with one block in \(\R^{d(f+1)}\) for every position and pair. The product
theorem for Ramsey sets~\cite[Theorem~20]{EGMRSS} and invariance under
dilation show that \(Z\) is Ramsey.
To each word \(\omega=(\omega_1,\ldots,\omega_N)\in(A^P)^N\)
associate the map \(\psi_\omega:B\to Z\) given by
\begin{equation}\label{eq:word}
 \psi_\omega(\mathbf{b})=\frac{1}{\sqrt{N|P|}}\,
   \bigl(\eta_{p'}(\mathbf{b},\omega_k(p'))\bigr)_
                  {(k,p')\in[N]\times P}.
\end{equation}
To compare \(\psi_\omega(\mathbf{b})\) and \(\psi_\omega(\mathbf{b}')\) for
\(\mathbf{b},\mathbf{b}'\in B\), fix a block \((k,p')\). Its point of \(A\)
is the same in both images, namely \(\omega_k(p')\). The isometry
in the first argument of \(\eta_{p'}\) gives squared distance
\(\norm{\mathbf{b}-\mathbf{b}'}^2\) before normalization. After normalization
each block contributes \(\norm{\mathbf{b}-\mathbf{b}'}^2/(N|P|)\).
There are \(N|P|\) blocks, so \(\psi_\omega\) is an isometry.

\smallskip\noindent\emph{A pair that works along a line.}
Take a coloring \(\chi:Z\to\Gamma\) with no rainbow copy of \(B\).
In their proof for hypercubes, Geh\'er, Sagdeev and
T\'oth~\cite[proof of Proposition~2]{GST} classify the event that a copy
is not rainbow by the pair of vertices that receive the same color.
They fix one endpoint and vary the other over a regular simplex to
bound the probability of this event. We use the pair itself as the
color of a word, so that the ordinary Hales--Jewett theorem will give
a whole line with the same pair.
Define the auxiliary coloring \(\widehat{\chi}:(A^P)^N\to P\).
For each word \(\omega\), let \(\widehat{\chi}(\omega)\) be the first pair
\(\{i,j\}\), in a fixed order on \(P\), for which the points
\(\psi_\omega(\mathbf{b}_i)\) and \(\psi_\omega(\mathbf{b}_j)\) have the same
\(\chi\)-color. Such a pair exists since \(\psi_\omega(B)\) is a copy
of \(B\). Take a monochromatic line for \(\widehat{\chi}\). Write
\(p=\{i,j\}\), \(i<j\), for its common color and \(I\) for its set
of variable positions.

Fix \(\mathbf{a}_0\in A\). For each \(\mathbf{a}\in A\), use on the line the
letter \(\theta_{\mathbf{a}}:P\to A\) that takes the value \(\mathbf{a}\) at
\(p\) and the value \(\mathbf{a}_0\) at every \(p'\ne p\). Denote the
resulting word by \(\omega_{\mathbf{a}}\in(A^P)^N\).
In~\eqref{eq:word}, the only blocks that can vary with \(\mathbf{a}\)
are \((k,p)\) with \(k\in I\). Even in these blocks the image
of \(\mathbf{b}_i\) is zero. Thus \(\psi_{\omega_{\mathbf{a}}}(\mathbf{b}_i)\)
is a single fixed point of \(Z\). Because every word on the line
has \(\widehat{\chi}\)-color \(p\), each \(\psi_{\omega_{\mathbf{a}}}(\mathbf{b}_j)\)
has the \(\chi\)-color of this fixed point.

It remains to identify the monochromatic set just obtained. Let
\(t=|I|\). For \(\mathbf{a},\mathbf{a}'\in A\), each of the \(t\) varying
blocks contributes the squared distance in~\eqref{eq:moving-vertex}
divided by \(N|P|\), and all other blocks contribute zero. Hence
\begin{equation}\label{eq:line-scale}
 \norm{\psi_{\omega_{\mathbf{a}}}(\mathbf{b}_j)-\psi_{\omega_{\mathbf{a}'}}(\mathbf{b}_j)}^2
 =\frac{t}{N|P|}\,\frac{N|P|}{q}\norm{\mathbf{a}-\mathbf{a}'}^2
 =\frac{t}{q}\norm{\mathbf{a}-\mathbf{a}'}^2 .
\end{equation}
Since \(1\leq t\leq N\) and \(t\) divides \(q\), the integer
\(h=q/t\) belongs to \([q]\). The monochromatic set is a copy of \(A_h\).
\end{proof}

\section{Returning to the original scale}\label{sec:scale}

\begin{proof}[Proof of \Cref{thm:main}]
If either \(A\) or \(B\) is a singleton, a one-point set suffices for
\(C\). Otherwise fix \(Z\) and \(q\) from \Cref{lem:scales}.
Geh\'er, Sagdeev and T\'oth~\cite[proof of Proposition~1]{GST}
color a layer by the position of a monochromatic edge to construct a
monochromatic prism. This product argument is used for cuboids
in~\cite{ShawCuboids}, for powers of a fixed simplex
in~\cite[Lemma~2.5]{GSXY}, and for two different factors
in~\cite[Lemma~7]{Shaw}. We give the induction with the monochromatic
factor allowed to vary among \(A_1,\ldots,A_q\).

\begin{Claim}\label{cl:layers}
For each positive integer \(M\) there is a finite configuration
\(C_M\) such that every coloring of \(C_M\) with no rainbow \(B\)
contains a monochromatic copy of \(A_{h_1}\times\cdots\times A_{h_M}\)
for some \(h_1,\ldots,h_M\in[q]\).
\end{Claim}

For \(M=1\) take \(C_1=Z\). Suppose \(C_M\) has been constructed.
As \(Z\) is Ramsey, the finite witness
principle~\cite[Proposition~4]{EGMRSS} gives a finite Euclidean set
\(D\) such that every coloring of \(D\)
with at most \(2^{|C_M|}\) colors has a monochromatic copy of \(Z\).
Set \(C_{M+1}=C_M\times D\) and fix a coloring
\(\chi:C_{M+1}\to\Gamma\) with no rainbow \(B\).

For each \(\mathbf{d}\in D\), the layer \(C_M\times\{\mathbf{d}\}\) is
isometric to \(C_M\) and contains no rainbow copy of \(B\).
Choose \(X_{\mathbf{d}}\subseteq C_M\) such that
\(X_{\mathbf{d}}\times\{\mathbf{d}\}\) is one of the monochromatic products
from the induction hypothesis. The map \(\mathbf{d}\mapsto X_{\mathbf{d}}\)
colors \(D\) by subsets of \(C_M\), so it uses at most \(2^{|C_M|}\)
colors. We obtain a copy \(Z'\subseteq D\) of \(Z\) on which
\(X_{\mathbf{d}}\) is the same set, say \(X\). Each layer
\(X\times\{\mathbf{d}\}\) with \(\mathbf{d}\in Z'\) is monochromatic,
although its \(\chi\)-color may depend on \(\mathbf{d}\).

Fix \(\mathbf{x}\in X\). In the copy \(\{\mathbf{x}\}\times Z'\) of \(Z\),
\Cref{lem:scales} gives \(Y\subseteq Z'\) congruent to some
\(A_h\) such that all points of \(\{\mathbf{x}\}\times Y\) have the same
\(\chi\)-color. For any \(\mathbf{x}'\in X\) and \(\mathbf{d}\in Y\),
the point \((\mathbf{x}',\mathbf{d})\) has the color of \((\mathbf{x},\mathbf{d})\),
since both belong to the same monochromatic layer. Thus \(X\times Y\)
is monochromatic. It is congruent to the product for \(X\) with one
more factor \(A_h\), proving the claim.

Ge, Shu, Xu and Yu~\cite[proof of Theorem~1.1, after Lemma~2.5]{GSXY}
recover a simplex at its original scale on the diagonal of a product of
identically scaled copies. To apply their observation to our finite
list of possible scales, take \(M=q^2\) and \(C=C_M\). In any coloring with no rainbow
\(B\), \Cref{cl:layers} gives an isometry from a product
\(A_{h_1}\times\cdots\times A_{h_M}\) onto a monochromatic subset of \(C\).
There are \(q^2\) indices \(h_1,\ldots,h_M\), each in \([q]\), so some
value \(h\) occurs at least \(q\) times, and hence at least \(h\) times.
Choose \(h\) such factors and fix a point in each remaining factor.
For \(\mathbf{a}\in A\), put \(\mathbf{a}/\sqrt h\) in every chosen factor.
For \(\mathbf{a},\mathbf{a}'\in A\), each of the \(h\) chosen factors contributes
\(\norm{\mathbf{a}-\mathbf{a}'}^2/h\) to the squared distance, and the fixed
factors contribute zero. The sum is \(\norm{\mathbf{a}-\mathbf{a}'}^2\).
The image of this diagonal under the isometry is the required
monochromatic copy of \(A\).
\end{proof}

\appendix

\section{A canonically Ramsey set that is not Ramsey}\label{appA}

Let \(T=\{0,1,2\}\subset\R\). Since no sphere contains three collinear
points, \(T\) is not Ramsey~\cite[Theorem~13]{EGMRSS}. Cheng and
Xu~\cite{CX} proved the canonical result for finite colorings of the
plane that are constant on each circle about the origin. We use
\Cref{thm:main} to force one monochromatic sphere. From this we obtain
six monochromatic spheres to which we can apply their observation
about consecutive squared radii.

\begin{Theorem}\label{thm:unit-progression}
There is a positive integer \(n\) such that \(\R^n\MR T\).
\end{Theorem}

\begin{proof}
Call a coloring \emph{bad} if every copy of \(T\) receives exactly two
colors. Suppose that \(\R^n\) admits a bad coloring for every \(n\).
A bad coloring contains no rainbow square of diagonal \(2\). Indeed,
the endpoints and the midpoint of each diagonal form a copy of \(T\).
If the square is rainbow, its center must have the color of an endpoint
of each diagonal, which is impossible.

Let \(V=\Q^{(\N)}\) be the space of finitely supported rational vectors
and put \(S_r=\{\mathbf{x}\in V:\norm{\mathbf{x}}^2=r\}\).
To obtain a bad coloring of \(V\) in which certain spheres are
monochromatic, it suffices to find such a coloring on each finite
\(F\subset V\), with the conditions restricted to \(F\).
Indeed, enumerate \(V\) and consider the partitions of each initial
segment that satisfy the conditions on that segment. These partitions
form a finitely branching tree under restriction. If every level is
nonempty, an infinite branch gives the required coloring.

We first make \(S_1\) monochromatic. Given a finite \(F\subset V\),
choose \(d\) so that \(F\subset\Q^d\), and enlarge \(F\cap S_1\) by
the points \(\pm\mathbf{e}_i\), \(i\in[d]\), where
\(\mathbf{e}_1,\ldots,\mathbf{e}_d\) are the standard basis vectors.
The resulting set \(A\) has full affine span
and is Ramsey by \Cref{thm:criterion}. Indeed, the coordinate functions
with weights \(1\) and the constant function \(1\) with weight \(-1\)
satisfy~\eqref{eq:certificate}, since all coordinates are rational and
their squares sum to \(1\) on \(A\). A square is also Ramsey, being a
Cartesian product of two Ramsey pairs~\cite[Theorem~20]{EGMRSS}.
Apply \Cref{thm:main} with \(A\) as the monochromatic target and a
square of diagonal \(2\) as the rainbow target. In a bad coloring of
\(\R^n\) for sufficiently large \(n\ge d\), this gives a monochromatic
copy of \(A\). Extend the isometry onto this copy to an isometric
embedding of \(\R^d\) into \(\R^n\) and pull the coloring back to \(F\).
By compactness, there is a bad coloring of \(V\) in which \(S_1\) is
monochromatic.

We record a geometric observation. Every positive rational \(a/b\),
where \(a,b\) are positive integers, is a sum of \(ab\) squares of
\(1/b\). Let \(r,s\) be positive rationals and \(t\in\Q\). Given
\(k\) mutually orthogonal vectors \(\mathbf{x}_1,\ldots,\mathbf{x}_k\in S_r\),
we can find \(\mathbf{y}\in S_s\) such that
\(\mathbf{x}_i\cdot\mathbf{y}=t\) for every \(i\), provided
\(kt^2\le rs\). Start with \((t/r)\sum_i\mathbf{x}_i\) and fill the
remaining squared norm in unused coordinates. In particular, if
\(2t^2\le rs\), any two points of \(S_r\) can be joined by an even
walk alternating between \(S_r\) and \(S_s\), with dot product \(t\)
on each edge. To see this, take a point of \(S_r\) orthogonal to both
given points. It has a common neighbor in \(S_s\) with each of them.

Suppose now that \(2\le r\le6\) is an integer and we have a bad
coloring of \(V\) in which \(S_1,\ldots,S_{r-1}\) are each
monochromatic. Write \(\alpha\) for the color of \(S_{r-1}\).
Two points of \(S_r\) at distance \(2\) have their midpoint in
\(S_{r-1}\). Hence at most one of them has color \(\alpha\), and if
neither does, they have the same color. The graph joining such pairs
in \(S_r\) is connected by the observation above with \(s=r\) and
\(t=r-2\), since \(2(r-2)^2\le r^2\) for \(2\le r\le6\).

Fix a finite \(F\subset V\) and choose a finite connected subgraph
with vertex set \(G\subset S_r\) containing \(F\cap S_r\). Take
\(m>|G|\). We construct linear isometries \(L_1,\ldots,L_m:V\to V\)
such that \(\norm{L_i\mathbf{x}-L_j\mathbf{x}}=2\) for
\(\mathbf{x}\in S_r\) and \(i\ne j\). Use \(r(r-2)\) common
coordinate blocks and \(2r\) further blocks for each \(i\in[m]\),
each block a copy of \(V\). Put \(\mathbf{x}/r\) in every common block
and in the \(2r\) blocks belonging to \(i\), and zero elsewhere.
The resulting vector \(L_i\mathbf{x}\) has squared norm
\(r^2\norm{\mathbf{x}/r}^2=\norm{\mathbf{x}}^2\). For \(i\ne j\)
and \(\mathbf{x}\in S_r\), the squared distance between
\(L_i\mathbf{x}\) and \(L_j\mathbf{x}\) is \(4r\norm{\mathbf{x}/r}^2=4\).

For each \(\mathbf{x}\in G\), at most one of its \(m\) images has color
\(\alpha\). Thus some \(L_iG\) avoids \(\alpha\) and is monochromatic
by connectivity. Pulling the coloring back along \(L_i\) makes
\(F\cap S_r\) monochromatic and preserves monochromaticity on each
earlier sphere, since \(L_i\) preserves norms. Compactness gives the
induction step. We may therefore fix a bad coloring \(\chi\) of \(V\)
in which \(S_1,\ldots,S_6\) are each monochromatic.

Following Cheng and Xu~\cite{CX}, we first determine the colors of
these six spheres. A point \(\mathbf{x}\in S_s\) is the midpoint of a
copy \(\{\mathbf{x}-\mathbf{u},\mathbf{x},\mathbf{x}+\mathbf{u}\}\)
of \(T\) whose endpoints have squared norms \(s+1\pm2t\), whenever
\(t\in\Q\) and \(t^2\le s\). Indeed, the geometric observation with
\(k=1\) gives a unit vector \(\mathbf{u}\) with
\(\mathbf{x}\cdot\mathbf{u}=t\). Taking \(s=j\) and \(t=0\) or
\(t=1/2\) shows that \(S_j\) differs in color from both \(S_{j+1}\)
and \(S_{j+2}\), whenever the indicated spheres occur among the six.
Now fix a point of \(S_{j+1/2}\), where \(j=1,2,3\). Taking
\(t=3/4\) and \(t=1/4\) gives two copies of \(T\) with this common
midpoint and with endpoints in \(S_j,S_{j+3}\) and
\(S_{j+1},S_{j+2}\), respectively. If the four spheres have distinct
colors, the midpoint color must belong to two disjoint pairs of colors.
It follows that \(S_j\) and \(S_{j+3}\) have the same color. Thus the
colors of \(S_1,\ldots,S_6\), in order, are \(a,b,c,a,b,c\), with
\(a,b,c\) distinct.

Put \(X=S_{4/3}\), \(Y=S_{5/3}\), and \(Z=S_{8/3}\). If
\(\norm{\mathbf{x}-\mathbf{y}}=1\), the four points
\(2\mathbf{x}-\mathbf{y},\mathbf{x},\mathbf{y},2\mathbf{y}-\mathbf{x}\)
form two consecutive copies of \(T\). Suppose the outer points have
distinct colors \(p,q\). If the two inner colors coincide, they differ
from both \(p\) and \(q\). If they differ, each outer color must be one
of the two inner colors. In either case,
\[
  \chi(\mathbf{y})=(p\,q)\chi(\mathbf{x}),
\]
where \((p\,q)\) interchanges \(p,q\) and fixes all other colors.
Since \(\norm{2\mathbf{x}-\mathbf{y}}^2
=2\norm{\mathbf{x}}^2-\norm{\mathbf{y}}^2+2\), the outer squared
norms are \(3,4\) for unit neighbors in \(X,Y\), and \(2,6\) for unit
neighbors in \(X,Z\). The corresponding transpositions are
\((c\,a)\) and \((b\,c)\).

Unit neighbors in \(X,Y\) have dot product \(1\). Since
\(2\le(4/3)(5/3)\), any two points of \(X\) are joined by an even path
of such edges. The transposition \((c\,a)\) cancels in pairs along the
path, so \(X\) has a single color \(\gamma\). Every point of \(Y\) or
\(Z\) has a unit neighbor in \(X\), by the observation with \(k=1\)
and dot product \(1\) or \(3/2\), respectively. Thus the colors of
\(Y\) and \(Z\) are \((c\,a)\gamma\) and \((b\,c)\gamma\).

If \(\gamma\notin\{a,b,c\}\), then \(Y\) and \(Z\) both have color
\(\gamma\). Taking \(s=5/3\) and \(t=0\) gives a monochromatic copy
of \(T\) with midpoint in \(Y\) and both endpoints in \(Z\).
If \(\gamma\in\{a,b,c\}\), take \(s=1\) and \(t=1/3\). This gives a
copy with midpoint in \(S_1\) and endpoints in \(X,Z\). Its colors
are \(\gamma,a,(b\,c)\gamma\). They coincide if \(\gamma=a\) and are
pairwise distinct if \(\gamma=b\) or \(c\). Both cases contradict the
choice of \(\chi\).
\end{proof}

\vspace{2.5mm}\
\noindent\textbf{Remark.}
Both \(T\) and \(\{0,1\}\) are canonically Ramsey, while their
Cartesian product is not. Indeed, consider the four-coloring
\(\mathbf{x}\mapsto\lfloor\norm{\mathbf{x}}^2\rfloor\bmod 4\)
from~\cite{CX}. If \(q_-,q_0,q_+\) are the squared norms of the points
of a copy of \(T\), in order, then \(q_-+q_+-2q_0=2\).
If all three colors agree, the corresponding combination of the
integer parts is a multiple of \(4\), while that of the fractional
parts lies strictly between \(-2\) and \(2\), a contradiction.
Every copy of \(T\times\{0,1\}\) therefore contains a nonmonochromatic
copy of \(T\), and its six points use at most four colors.

\vspace{2.5mm}\
\noindent\textbf{Remark.}
Similar arguments can be used to show that $\{0,t,1\}\subset\R$ is canonically Ramsey for every $t \in \Q$. The case of $t\not\in\Q$ seems different: one can construct, for every $n$, a coloring of $(\Q(t))^n$ containing neither monochromatic nor rainbow copies of $\{0,t,1\}$. This however is not sufficient to conclude that the set is not canonically Ramsey.

\end{document}